\documentclass[11pt]{article}
\usepackage{amsmath,amssymb,amscd,latexsym,amsthm,mathrsfs}
\usepackage{hypbmsec,longtable}
\usepackage{graphicx}
\usepackage{tikz}\usetikzlibrary{positioning}
\usepackage[a4paper,left=2.25cm,right=2.25cm,top=3cm,bottom=3cm,headsep=1cm]{geometry}
\usepackage[unicode]{hyperref}
\makeatletter
\newcommand{\gettikzxy}[3]{
  \tikz@scan@one@point\pgfutil@firstofone#1\relax
\pgfmathsetmacro{#2}{\the\pgf@x/\linkpatternunit}
\pgfmathsetmacro{#3}{\the\pgf@y/\linkpatternunit}
}
\tikzset{label anchor/.code={%
    \let\tikz@auto@anchor=\pgfutil@empty
    \def\tikz@anchor{#1}
  },
  label anchor/.default=center
}
\makeatother
\tikzset{arrow/.style={postaction={decorate,thick,decoration={markings,mark = at position #1 with {\arrow{>}}}}},arrow/.default=0.5}
\tikzset{invarrow/.style={postaction={decorate,thick,decoration={markings,mark = at position #1 with {\arrow{<}}}}},invarrow/.default=0.5}
\newdimen\linkpatternunit%
\newcount\linkpatternsize%
\newcount\lpsize
\newif\iflinkpatterninverted
\newif\iflinkpatterntikzstarted
\newif\iflinkpatternboxed
\newif\iflinkpatternaxis
\newif\iflinkpatternstraightlines
\newif\iflinkpatternnumbered
\newif\iflinkpatternalias
\newif\iflinkpatternnode
\newif\iflinkpatterncentered
\newcount\linkpatternfused
\pgfkeys{/linkpattern/.cd,centered/.is if=linkpatterncentered,inverted/.is if=linkpatterninverted,numbered/.is if=linkpatternnumbered,tikzstarted/.is if=linkpatterntikzstarted,straight lines/.is if=linkpatternstraightlines,boxed/.is if=linkpatternboxed,axis/.is if=linkpatternaxis,vertexcolor/.store in=\linkpatternvertexcolor,edgecolor/.store in=\linkpatternedgecolor,boxcolor/.code={\linkpatternboxedtrue\def\linkpatternboxcolor{#1}},tikzoptions/.style={every linkpattern/.append style={#1}},size/.code={\linkpatternsize=#1},numbering0/.code={\def\lpnumbering{#1}\def\linkpatternnumbering{#1}},numbering/.style={numbered,numbering0={#1}},unit/.code={\linkpatternunit=#1},height/.store in=\linkpatternheight,shape/.store in=\linkpatternshape,looseness/.store in=\linkpatternlooseness,squareness/.store in=\linkpatternsquareness,extra space/.store in=\linkpatternextraspace,width/.store in=\linkpatternwidth,alias/.is if=linkpatternalias,pos/.store in=\linkpatternpos,labeloptions/.style={labeloptionslist/.append style={#1}},labeloptionslist/.style={inner sep=2pt,font=\scriptsize,execute at begin node=$,execute at end node=$,label anchor={#1+180}},nodeon/.is if=linkpatternnode,node/.style={nodeon,nodeoptionslist/.append style={#1}},nodeoptionslist/.style={},
pipedream/.style={shape=pipedream,looseness=0,straight lines,numbering0=tangle},
tangle/.style={shape=tangle,numbering0=tangle},
every linkpattern/.style={x=\linkpatternunit,y=\linkpatternunit},
fused/.code={\linkpatternfused=#1},
vertex/.style={circle,thin,solid,draw=black,fill=\linkpatternvertexcolor,inner sep=1.5pt,draw opacity=1,transform shape},
edge/.style={very thick,solid,draw=\linkpatternedgecolor,draw opacity=1}}
\linkpatterncenteredfalse%
\linkpatterninvertedfalse%
\linkpatternnumberedfalse%
\linkpatterntikzstartedfalse%
\linkpatternboxedfalse%
\linkpatternaxistrue%
\linkpatternaliastrue%
\linkpatternstraightlinesfalse%
\def\linkpatternlooseness{0.2}
\def\linkpatternsquareness{0.35}
\def\linkpatternvertexcolor{red}%
\def\linkpatternedgecolor{blue}%
\def\linkpatternboxcolor{none}%
\def\linkpatternheight{0}
\def\linkpatternwidth{0}
\def\linkpatternshape{default}
\def\linkpatternnumbering{default}
\def\linkpatternpos{(0,0)}
\def\linkpatternextraspace{0}
\linkpatternnodefalse
\def\firstchar#1#2\empty{#1}%
\def\linkpatterndo#1#2{
\edef\param{\csname linkpattern#2\endcsname}
\edef\firstcharparam{\expandafter\firstchar\param\empty}
\expandafter\ifcat\firstcharparam a
\expandafter\ifx\csname linkpattern#1\param\endcsname\relax
\csname linkpattern#1unknown\endcsname
\else
\csname linkpattern#1\csname linkpattern#2\endcsname\endcsname
\fi
\else
\csname linkpattern#1unknown\endcsname
\fi
}%
\def\linkpatterncoordtangle{\ifnum\x>\lphalfsize\pgfmathparse{\lpsize+1-\x}\xdef\lpcoordx{\pgfmathresult}\xdef\lpcoordy{\lpheight}\xdef\lpangle{270}\else\xdef\lpcoordx{\x}\xdef\lpcoordy{-\lpheight}\xdef\lpangle{90}\fi}
\def\linkpatterncoordpipedream{\ifnum\x>\lphalfsize\pgfmathparse{\lpsize+1-\x-0.5}\xdef\lpcoordx{\pgfmathresult}\xdef\lpcoordy{0}\xdef\lpangle{270}\else\pgfmathparse{0.5-\x}\xdef\lpcoordy{\pgfmathresult}\xdef\lpcoordx{0}\xdef\lpangle{0}\fi}
\def\linkpatterncoordrectangle{
\ifnum\x>\lptqsize
\pgfmathparse{\lpsize+1-\x-0.5}\xdef\lpcoordx{\pgfmathresult}\xdef\lpcoordy{0}\xdef\lpangle{270}
\else\ifnum\x>\lphalfsize
\pgfmathparse{\x-\lptqsize-0.5}\xdef\lpcoordy{\pgfmathresult}\xdef\lpcoordx{\linkpatternwidth}\xdef\lpangle{180}
\else\ifnum\x>\linkpatternheight
\pgfmathparse{\x-\linkpatternheight-0.5}\xdef\lpcoordx{\pgfmathresult}\xdef\lpcoordy{-\linkpatternheight}\xdef\lpangle{90}
\else
\pgfmathparse{0.5-\x}\xdef\lpcoordy{\pgfmathresult}\xdef\lpcoordx{0}\xdef\lpangle{0}
\fi\fi\fi
}%
\def\linkpatternsetsizeunknown{
\global\lpsize=\linkpatternsize
\if\linkpatternheight0
\xdef\maxsep{0}
\foreach \x/\xx in \mylist%
{%
\edef\tempx{\withoutprime{\x}}
\edef\tempxx{\withoutprime{\xx}}
\pgfmathparse{max(\maxsep,abs(\tempx-\tempxx))}
\xdef\maxsep{\pgfmathresult}
}%
\pgfmathparse{0.25+0.8*\linkpatternsquareness*\maxsep}
\xdef\lpheight{\pgfmathresult}
\else
\xdef\lpheight{\linkpatternheight}
\fi
}
\newcount\tempsize
\def\linkpatternrightmostunknown{
\global\lpsize=0
\global\tempsize=0
\foreach\x/\labx in \linkpatternnumbering
{
\edef\tempx{\withoutprime{\x}}
\ifnum\lpsize<\tempx\global\lpsize=\tempx\fi
\global\advance\tempsize by 1
}
\ifnum\tempsize>\lpsize\global\lpsize=\tempsize\fi
}%
\def\linkpatternrightmostdefault{
\global\lpsize=0
\global\tempsize=0
\foreach \x/\y in \mylist
{
\edef\tempx{\withoutprime{\x}}
\ifnum\lpsize<\tempx\global\lpsize=\tempx\fi
\ifx\x\y
\global\advance\tempsize by 1
\else
\edef\tempy{\withoutprime{\y}}
\ifnum\lpsize<\tempy\global\lpsize=\tempy\fi%
\global\advance\tempsize by 2
\fi
}
\ifnum\tempsize>\lpsize\global\lpsize=\tempsize\fi
}%
\def\linkpatternrightmosttangle{
\global\lpsize=0
\global\tempsize=0
\foreach \x/\y in \mylist
{
\edef\tempx{\withoutprime{\x}}
\ifnum\lpsize<\tempx\global\lpsize=\tempx\fi
\ifx\x\y
\global\advance\tempsize by 1
\else
\edef\tempy{\withoutprime{\y}}
\ifnum\lpsize<\tempy\global\lpsize=\tempy\fi%
\global\advance\tempsize by 2
\fi
}
\global\advance\lpsize by\lpsize
\ifnum\tempsize>\lpsize\global\lpsize=\tempsize\fi
}%

\newcommand\linkpattern[2][]{
{
\pgfkeys{/linkpattern/.cd,#1}
\edef\mylist{#2}
\def\primetest##1'{}%
\def\hasaprime##1{\expandafter\primetest##1''}
\def\internalwithoutprime##1'{##1}%
\def\withoutprime##1{\if\hasaprime##1 %
\expandafter\internalwithoutprime##1\else ##1\fi}%
\iflinkpatternnumbered%
\iflinkpatterninverted
\tikzset{/linkpattern/lbl/.style n args={3}{label={[/linkpattern/labeloptionslist=-##1,##3] ##1:##2}}}%
\else%
\tikzset{/linkpattern/lbl/.style n args={3}{label={[/linkpattern/labeloptionslist=##1,##3] ##1:##2}}}%
\fi%
\else%
\tikzset{/linkpattern/lbl/.style={}}%
\fi%
\tikzifinpicture{\linkpatterntikzstartedtrue%
\begin{scope}[shift=\linkpatternpos,/linkpattern/every linkpattern]
}{%
\linkpatterntikzstartedfalse%
\iflinkpatterncentered
\begin{tikzpicture}[baseline=(current  bounding  box.center),/linkpattern/every linkpattern]%
\else
\begin{tikzpicture}[baseline=0,/linkpattern/every linkpattern]%
\fi
}%
\begin{scope}[local bounding box=link pattern box]
\iflinkpatterninverted%
\begin{scope}[yscale=-1]%
\fi%
\linkpatterndo{setsize}{shape}
\ifnum\lpsize=0
\linkpatterndo{rightmost}{numbering}
\fi
\pgfmathtruncatemacro{\lphalfsize}{\lpsize/2}
\linkpatterndo{numbering}{numbering}
\iflinkpatternboxed
\linkpatterndo{drawbox}{shape}
\else
\iflinkpatternaxis
\linkpatterndo{drawaxis}{shape}
\fi
\fi
\foreach\xx/\xlab/\opt in \lpnumbering
{
\ifx\xlab\opt\def\opt{}\fi
\if\hasaprime\xx %
\pgfmathtruncatemacro{\xx}{\lpsize+1-\withoutprime{\xx}}
\fi
\ifnum\linkpatternfused>1
\pgfmathsetmacro{\x}{0.4*(0.5+\linkpatternfused*(0.5+floor((\xx-1)/\linkpatternfused)))+0.6*\xx}
\else
\def\x{\xx}
\fi
\linkpatterndo{coord}{shape}
\iflinkpatternalias\def\xlabb{\xlab}\else\def\xlabb{\xx}\fi
\path (\lpcoordx,\lpcoordy) coordinate[/linkpattern/vertex,/linkpattern/lbl={\lpangle+180}{\xlab}{\opt},alias=v\xlabb] (v\xx) ++(\lpangle:\linkpatternunit) coordinate[alias=vv\xlabb] (vv\xx); 
}
\foreach \a/\b/\c in \mylist
{
\if\hasaprime\a %
\pgfmathtruncatemacro{\a}{\lpsize+1-\withoutprime{\a}}
\fi
\ifx\b\c\def\c{}\fi
\draw[/linkpattern/edge]
\ifx\a\b
(v\a)
\c
--
++(0,\lpheight);
\else
\pgfextra{
\if\hasaprime\b %
\pgfmathtruncatemacro{\b}{\lpsize+1-\withoutprime{\b}}
\fi
\gettikzxy{(v\a)}{\ax}{\ay}
\gettikzxy{(v\b)}{\bx}{\by}
\gettikzxy{(vv\a)}{\axx}{\ayy}
\gettikzxy{(vv\b)}{\bxx}{\byy}
\pgfmathsetmacro{\dist}{sqrt((\ax-\bx)*(\ax-\bx)+(\ay-\by)*(\ay-\by))}
\pgfmathsetmacro{\abx}{(\axx-\ax)*\dist*\linkpatternsquareness+(\bx-\ax)*\linkpatternlooseness)}
\pgfmathsetmacro{\aby}{(\ayy-\ay)*\dist*\linkpatternsquareness+(\by-\ay)*\linkpatternlooseness)}
\pgfmathsetmacro{\bax}{(\bxx-\bx)*\dist*\linkpatternsquareness+(\ax-\bx)*\linkpatternlooseness)}
\pgfmathsetmacro{\bay}{(\byy-\by)*\dist*\linkpatternsquareness+(\ay-\by)*\linkpatternlooseness)}
}
(v\a)
\c
\iflinkpatternstraightlines
\pgfextra{
\pgfmathsetmacro{\t}{((\ax-\bx)*\bay-(\ay-\by)*\bax)/(\aby*\bax-\abx*\bay)}
\pgfmathsetmacro{\abx}{\t*\abx}
\pgfmathsetmacro{\aby}{\t*\aby}
}
[rounded corners=0.2\linkpatternunit] -- ++(\abx,\aby) -- (v\b);
\else
.. controls ++(\abx,\aby) and ++(\bax,\bay) .. 
\fi
(v\b);
\fi
}
\end{scope}
\iflinkpatternnode
\node[fit=(link pattern box),/linkpattern/nodeoptionslist] {};
\fi
\iflinkpatterninverted
\end{scope}
\fi
\iflinkpatterntikzstarted
\end{scope}
\else%
\end{tikzpicture}%
\fi%
}}%
\newcommand\tanglelinkpattern[3][]{%
{
\pgfkeys{/linkpattern/.cd,#1}
\iflinkpatterninverted
\begin{tikzpicture}[/linkpattern/every linkpattern,baseline=\linkpatternunit]%
\else
\begin{tikzpicture}[/linkpattern/every linkpattern,baseline=-\linkpatternunit]%
\fi
\linkpattern[#1,tikzstarted,numbered=false]{#3}
\pgfmathtruncatemacro{\lptempsize}{2*\linkpatternsize}
\iflinkpatterninverted
\begin{scope}[yshift=0.5*\linkpatternunit]
\else
\begin{scope}[yshift=-0.5*\linkpatternunit]
\fi
\linkpattern[tangle,#1,tikzstarted,size=\lptempsize,
numbering=halftangle,
height=0.5]{#2}
\end{scope}
\end{tikzpicture}%
}}
\newcommand\diag[4][]{%
\pgfkeys{/linkpattern/.cd,#1}
\iflinkpatterntikzstarted\else%
\begin{tikzpicture}[scale=0.5]
\fi%
\iflinkpatterninverted%
\begin{scope}[yscale=-1]%
\fi%
\draw (0,0) grid (#2,#3);
\edef\mylist{#4}
\foreach\y/\x/\z in \mylist
{
\ifx\x\z
\draw[decorate,decoration={zigzag,
amplitude=1pt,segment length=5pt}]
(\x-0.5,#3) -- (\x-0.5,\y-0.5) node[circle,fill=black,inner sep=2pt] {} -- (#2,\y-0.5);
\else
\node at (\x-0.5,\y-0.5) {$\z$};
\fi
}
\iflinkpatterninverted
\end{scope}
\fi
\iflinkpatterntikzstarted\else%
\end{tikzpicture}%
\fi%
}
\makeatletter
\tikzset{circle split part fill/.style  args={#1,#2}{%
 alias=tmp@name,
  postaction={%
    insert path={
     \pgfextra{%
     \pgfpointdiff{\pgfpointanchor{\pgf@node@name}{center}}%
                  {\pgfpointanchor{\pgf@node@name}{east}}%
     \pgfmathsetmacro\insiderad{\pgf@x}
      \fill[#1] (\pgf@node@name.base) ([xshift=-\pgflinewidth]\pgf@node@name.east) arc
                          (0:180:\insiderad-\pgflinewidth)--cycle;
      \fill[#2] (\pgf@node@name.base) ([xshift=\pgflinewidth]\pgf@node@name.west)  arc
                           (180:360:\insiderad-\pgflinewidth)--cycle;                    }}}}}  
 \makeatother
\tikzset{bdot/.style={circle,circle split,draw,circle split part fill={black,white},thin,inner sep=1pt}}%
\tikzset{wdot/.style={circle,circle split,draw,circle split part fill={white,black},thin,inner sep=1pt}}%
\newcommand\circlelinkpattern[2][]{
{
\pgfkeys{/linkpattern/.cd,#1}
\iflinkpatterntikzstarted\else%
\begin{tikzpicture}[/linkpattern/every linkpattern]%
\fi%
\iflinkpatterninverted%
\begin{scope}[yscale=-1]%
\fi%
\global\lpsize=\linkpatternsize
\edef\mylist{#2}
\foreach \x/\y in \mylist
{
\ifnum\x>\lpsize\global\lpsize=\x\fi
\ifnum\y>\lpsize\global\lpsize=\y\fi
}
\iflinkpatternaxis
\draw (0,0) circle (1);
\fi
\foreach\x in {1,...,\lpsize}
{
\pgfmathparse{(0.3*floor((\x-1)/\linkpatternfused)+0.7*((\x-0.5)/\linkpatternfused-0.5))*\linkpatternfused*360/\lpsize}
\coordinate[/linkpattern/vertex] (v\x) at (\pgfmathresult:1);
}
\foreach \x/\y/\z in \mylist
{
\ifx\y\z%
\draw[/linkpattern/edge] (v\x) .. controls ($0.5*(v\x)$) and  ($0.5*(v\y)$) .. (v\y);
\else
\draw[/linkpattern/edge] \z (v\x) .. controls ($0.5*(v\x)$) and  ($0.5*(v\y)$) .. (v\y);
\fi
}
\iflinkpatternnumbered%
\pgfmathparse{\lpsize/\linkpatternfused}
\global\lpsize=\pgfmathresult
\def\linkpatternnumbering{1,...,\lpsize}
\newdimen\angle
\foreach\x/\xx/\opt in \linkpatternnumbering
{
  \pgfmathsetmacro{\angle}{360/\lpsize*(\x-1)}
\ifx\xx\opt%
  \node[outer sep=1pt,anchor=180+\angle] at (\angle:1) {$\scriptstyle\xx$}; 
\else
  \node[outer sep=1pt,anchor=180+\angle,\opt] at (\angle:1) {$\scriptstyle\xx$}; 
\fi
}
\fi%
\iflinkpatterninverted%
\end{scope}
\fi%
\iflinkpatterntikzstarted\else%
\end{tikzpicture}%
\fi%
}}%
\newdimen{\loopcellsize}
\tikzset{bgplaq/.style={draw=black,fill=\linkpatternboxcolor}}
\def\plaqwest{}
\def\plaqeast{}
\def\plaqnorth{}
\def\plaqsouth{}
\pgfkeys{/linkpattern/.cd,west/.store in=\plaqwest,east/.store in=\plaqeast,north/.store in=\plaqnorth,south/.store in=\plaqsouth}
\def\plaqname{plaq}
\newcommand\plaq[2][]{
\node[bgplaq,rectangle,draw,use as bounding box,minimum size=\loopcellsize,transform shape] (\plaqname) {};
\pgfkeys{/linkpattern/.cd,#1}
\ifx#2\empty\else
\begin{scope}[x=\loopcellsize,y=\loopcellsize]
\csname plaq#2\endcsname
\end{scope}\fi
}
\tikzset{loop/.code={\def\plaqname{loop-\the\pgfmatrixcurrentrow-\the\pgfmatrixcurrentcolumn}},loop/.append style={matrix,row sep={\loopcellsize,between origins},column sep={\loopcellsize,between origins}}}
\linkpatterninvertedtrue
\ifx\pdfoutput\undefined\else
\hypersetup{
colorlinks=true,
linkcolor=blue,
citecolor=blue,
urlcolor=blue,
filecolor=blue,
bookmarksnumbered=true,
pdfstartview=FitH,
pdfhighlight=/N
}
\fi

\theoremstyle{remark}
\newtheorem*{remark*}{Remark}

\newcommand{\BE}{\begin{equation}}
\newcommand{\EE}{\end{equation}}

\renewcommand{\author}[1]{\large\rm #1\\ \bigskip}
\newcommand{\address}[1]{{\normalsize\it #1\\}\bigskip}
\renewcommand{\title}[1]{\bigskip\bigskip\Large\bf #1\bigskip\bigskip\\}

\begin{document}

\vglue .3 cm

\begin{center}

\title{$1324$-avoiding permutations revisited}
\author{\renewcommand{\thefootnote}{\fnsymbol{footnote}}
Andrew R. Conway\footnote{email:  {\tt andrew1324@greatcactus.org}  }    Anthony J. Guttmann\footnote{email: {\tt tonyg@ms.unimelb.edu.au}} and Paul Zinn-Justin\footnote{email:  {\tt pzinn@unimelb.edu.au}  } }

\address{ 
School of Mathematics and Statistics,\\
The University of Melbourne, Victoria 3010, Australia}

\end{center}
\setcounter{footnote}{0}

\begin{abstract}
We give an improved algorithm for counting the number of $1324$-avoiding permutations, resulting in $14$ further terms of the generating function, which is now known for all lengths $\le 50$. We re-analyse the generating function and find additional evidence for our earlier conclusion that unlike other classical length-$4$ pattern-avoiding permutations, the generating function does not have a simple power-law singularity, but rather, the number of $1324$-avoiding permutations of length $n$  behaves as 
\[
B\cdot \mu^n \cdot \mu_1^{\sqrt{n}} \cdot n^g.
\]
We estimate $\mu=11.600 \pm 0.003$, $\mu_1 = 0.0400 \pm 0.0005$, $g = -1.1 \pm 0.1$ while the estimate of $B$ depends sensitively on the precise value of $\mu$, $\mu_1$ and $g$. This reanalysis provides substantially more compelling arguments for the presence of the stretched exponential term $\mu_1^{\sqrt{n}}$.

\end{abstract}

\section{Introduction}
In an earlier paper \cite{CG15}, two of the current authors gave further coefficients and a detailed analysis of the generating function for 1324 pattern-avoiding permutations (PAPs), extending the known ordinary generating function (OGF) by a further 5 terms. That analysis led us to conjecture that, unlike the known length-4 PAPs, notably the classes $Av(1234)$ \cite{IG90} and $Av(1342)$ \cite{MB97},  the OGF for $Av(1324)$ included a stretched exponential term. That is to say, if $p_n$ denotes the number of $n$-step $Av(1324)$ permutations, then 
\begin{equation} \label{eqn:pn}
p_n \sim B\cdot \mu^n \cdot \mu_1^{\sqrt{n}} \cdot n^g,
\end{equation}
 where estimates of the parameters were given.

In the present paper, we present a new, substantially improved algorithm that allows us to give 14 further terms\footnote{The only limit to obtaining additional terms is computer memory. The present calculation required about 4.2TB of (distributed) memory.}.

This stretched exponential behaviour is not without precedent. 
There are a number of models in mathematical physics whose coefficients possess this more complex asymptotic structure. In particular, Duplantier and Saleur \cite{DS87} and Duplantier and David \cite{DD88} studied the case of {\em dense}\/ polymers in two dimensions, and found the partition functions had the asymptotic form $const \cdot \mu^n \cdot \mu_1^{n^\sigma} \cdot n^g.$ In \cite{OPB93}, Owczarek, Prellberg and Brak investigated an exactly solvable model of interacting partially-directed self-avoiding walks (IPDSAW), and found the coefficients behaved with this asymptotic form, and estimated $\sigma = 1/2$, $g = -3/4$, while the sub-exponential growth constant $\mu_1$ was found to more than 5 digit accuracy. From \cite{BGW92} the value of $\mu$ is exactly known. Subsequently Duplantier \cite{D93} pointed out that $\sigma = 1/2$ is to be expected, not only for IPDSAWs, but also for SAWs in the collapsed regime. He went on to predict the value of the exponent $g$ in that case. For self-avoiding walks and polygons attached to a surface and pushed toward the surface by a force applied at their top vertex, Beaton {\em et al}\/ \cite{BGJL15} gave probabilistic arguments for stretched exponential behaviour, but with growth $\mu_1^{n^{3/7}}$.

Such stretched exponential behaviour is also seen in other combinatorial problems. If one considers the cogrowth series of certain infinite, finitely generated amenable groups \cite{EPG17}, one sees similar, and sometimes more complex, behaviour. For example, for the lamplighter group, the coefficients of the cogrowth series $l_n$ behave as 
\[
l_n \sim const.\cdot 9^n \cdot \mu_1^{n^{1/3}} \cdot n^{1/6},
\]
\cite{R03}, whereas for the wreath products ${\mathbb Z} \wr {\mathbb Z}$ and $({\mathbb Z} \wr {\mathbb Z}) \wr {\mathbb Z}$ the coefficients behave as 
\[
w_n \sim const.\cdot 16^n \cdot \mu_1^{n^{1/3}\log^{2/3}(n)} \cdot n^{g},
\] 
and  
\[
w_n \sim const.\cdot 36^n \cdot \mu_1^{\sqrt{n\log(n)}} \cdot n^{g},
\]
respectively \cite{PS-C02}. For the case of $1324$-PAPs we have good numerical evidence for the absence of such a confluent logarithmic term in the exponent, which we discuss below.

In the next section we give details of the algorithm. In subsequent sections we give our analysis.

\section{The algorithm}
The algorithm like many is based upon recursive solution of a set of equations
\[
f(S)=\sum_{s \in n(S)} f(s)
\]
where $n(S)$ is a set (or possibly multiset if the same $s$ appears with multiplicity) 
of possible substates of $S$, culminating
in some final states for which $f(S)=1$. These states correspond to the build up of permutations
one entry at a time, with each pass through the equation corresponding to adding one extra entry. 

As an example, one could use this formalism to enumerate all $n$ length permutations, saying the state is the number $n$ of
as yet unchosen elements. Then state $n$ would have $n$ substates, each $n-1$. This reduces to the normal factorial recurrence $f(n)=n f(n-1)$.
To enumerate PAPs, a more complex state is needed.

In the prior paper \cite{CG15}, the state consisted of a series of numbers being the length of contiguous series of
unchosen elements of the permutation, together with brackets to store sufficient information to prevent
a $1324$ pattern. For instance, the state $4(2)1$ means that there are four contiguous numbers left to be chosen,
then two contiguous numbers that may not be chosen until all numbers after them have been chosen, then 1 number.
Each pass through the equation reduced the sum of the available numbers by one.

The big insight is that this is unnecessarily fine grained. If one adds together all states with a given pattern
of brackets and numbers, ignoring what the numbers actually are other than their total sum, the equations still
work. As the start state and end state contain only a single number, there is no problem in summing them.
This significantly reduces the total number of states, 
and thus the running time and memory use of the enumeration.

A second insight is that by tracking states by what is taken rather than by what is available,
the states can be represented by Dyck paths, or link patterns, or by
any set in bijection with these, and enumerated by the Catalan numbers \cite{Stanley}; 
here we choose to use link patterns, as they provide a convenient graphical
description of the algorithm.
Since explicit bijections of these objects in length say $2k$ 
to $\{1,\ldots,\text{Catalan}_k=\frac{(2k)!}{k!(k+1)!}\}$ are known, we can encode them as integers.
The reader should be warned that $k$ here is {\em not}\/ the length of the permutation; it can vary
from state to state, with the upper bound $2k\le n$.

The intuitive
motivation of this way of tracking states comes from considering a prefix $P$ of a $1324$ avoiding permutation of $1 \dots n$, and
considering what constraints it puts on subsequent elements of the permutation. If $P$ contains $n$, then this cannot be part of a $1324$ pattern
where the $4$ is in the suffix, so the $n$ is irrelevant as a constraint on the future, and can be safely ignored, turning the 
problem into a permutation of $1 \dots n-1$. This process is repeated 
until all numbers in $P$ are lower than the largest number remaining in the suffix.
The remaining numbers in $P$ must be $132$ avoiding, as otherwise the largest number, now in the suffix, would cause
a $1324$ pattern. $132$ PAPs of a given length are enumerated by the Catalan numbers, and are readily bijectable to link patterns.

\subsection{States}\label{sec:state}
The algorithm proceeds in time steps, $t=0,\ldots,n$, where $n$ is the length of the $1324$ permutations
to count. We can also parameterise steps by the number $s$ of elements left to insert, with $s=n-t$.
This number will decrease by one each step through the algorithm, 
and will generally be assumed to be present hereafter and not
explicitly mentioned.

A state of the algorithm consists of a link pattern,
that is a {\em matching}\/ of $2k$ vertices on a line (i.e., a partition of these vertices into pairs) 
in such a way that the links (the pairings) can be drawn in a half-plane without crossing each other.
The number $k$ is state-dependent (we shall give some bounds on $k$ below).

Each link in the link pattern represents a consecutive (but not necessarily ordered) 
sequence of elements of the permutation. Their actual composition, length and order, is irrelevent to 
future constraints on what can be inserted where, as long as no $132$ sequences were 
created during its construction. There are two things that matter about it:
\begin{itemize}
	\item Its location relative to other gaps and consecutive elements. This is recorded by the right end of the link.
	\item The location of the smallest element inserted prior to one of these elements. This is recorded by the left end of the link.	      
\end{itemize}
The link then represents a potential $13$ portion of the $1324$ pattern. In the future, no insertions will be allowed
inside this link (a potential $2$) until everything to the right of this link (a potential $4$) is complete,
at which point this link will be discarded as no longer relevent.

New elements in the permutation can be added in the implicit gaps either at the start of the pattern or after each link, unless as just mentioned,
some other link surrounds it. This would mean that the element would be a $2$ in the pattern, and then a $4$ would inevitably arise as the state
always assumes the largest element has not arrived (by simplifying the problem as described earlier if it has arrived).

As an example of the state from a concrete permutation, see table~\ref{table:eg1fullpermutation}.

\begin{table}
	\begin{tabular}{lp{6.5cm}l}
		Element  & Notes & Result \\ \hline
		& Start state & $\varnothing$ \\
		$5$ & Not consecutive with anything; future elements could go either side. (1)& $5\atop\linkpattern{1/2}$ \\
		$4$ & Consecutive with 5 and so merged into it. (3) & $4-5\atop\linkpattern{1/2}$ \\
		$2$ & New link as not consecutive with anything. (1) & $~~2~~~~4-5\atop\linkpattern{1/2,3/4}$ \\
		$7$ & Larger than a previous link, makes constraint that no new elements between $2$ and $7$ may be added until every element greater than $7$ has been added. (1)& $~~~~~~2~~~4-5~~7\atop\linkpattern{1/6,2/3,4/5}$ \\
		$10$ & Removed from consideration as largest element. (3) & $~~~~~~2~~~4-5~~7\atop\linkpattern{1/6,2/3,4/5}$ \\
		$8$ & Merged with $7$. (2)  & $~~~~2~~~~4-5~7-8\!\!\!\!\!\atop\linkpattern{1/6,2/3,4/5}$ \\
		$9$ & Merges the $7-8$ link with the largest element; said link removed from consideration. (4) & $~~2~~~~4-5\atop\linkpattern{1/2,3/4}$ \\
		$1$ & Merges with $2$ link. (3) & $1-2~~4-5\atop\linkpattern{1/2,3/4}$ \\
		$3$ & Merges the $1-2$ and $3-5$ links. (4) & $1-5\atop\linkpattern{1/2}$ \\
		$6$ & Merges the $1-5$ link with the largest element. (4) & $\varnothing$ \\
	\end{tabular}
	\caption{An example of the state as the permutation $5~4~2~7~10~8~9~1~3~6$ is built up. The numbers above the link diagrams indicate the actual numbers that the link end represents. The numbers in parentheses correspond to the four types of insertion given in the text.
	   Conversely, consider the similar permutation $5~4~2~7~10~6~9~1~3~8$ which is not $1324$ avoiding. The first five elements are the same; the sixth element is not allowed, 
	   as it would would have to go inside the loop ending at $7$.}
	\label{table:eg1fullpermutation}
\end{table}

Generally, to get all possible child states of a given state one could insert elements into any
of the spaces before or after links that are not surrounded by other links. These elements may or may not
be consecutive with the links on the left and/or right, giving in general four possibilities for each insertion point: 
\begin{enumerate}
	\item Insert an element not consecutive on either side. This will add a link from the inserted place back to the start of the pattern. This increases the number of links, 
	      and hence implicit gaps, by one. This represents that a gap has been split into two gaps.
	\item Insert an element consecutive with the link to the left. This may not be done if inserting to the left of the pattern. 
	      The new element merges in with the link to its immediate left, so does not change its right end. However, it is possible that this changes the
	      smallest element inserted before this link. To reflect this, extend the left end of this link back to the start of the pattern.
	\item Insert an element consecutive with a link opening to the right (or the end of the pattern). 
	      This will add a constraint from the start of the pattern to the elements represented by the link opening it is connected to. 
	      Represent this by moving said opening to the start of the pattern. 
	      Note that a series of consecutive link openings all represent the same point, so all consecutive link openings immediately to the
	      right of the insertion point would have their openings moved to the start of the pattern. 
	      One may also insert an element contiguous with the right of the rightmost gap, in which case this is the highest number and has no effect on the link pattern.
	\item Insert an element consecutive with a link to both left and right (the right may also be the end of the pattern). This may not be done if inserting to the left of the pattern. 
	      This will merge the links on either side of the inserted element. As in the prior option, for all consecutive link openings immediately to the right of the insertion point, move their openings to the start of the pattern. The constraint implied by the link to the left is now covered by the links to the right, and the gap to the right of said link no longer exists. So erase the link to the left. This reduces the number of links, and hence implicit gaps, by one, which is appropriate as this element being consecutive with both sides of the
	      gap means that the gap has been filled.
	      If there are no links to the right, then this just fills in the last space, and erasing the link to the left has dealt with everything. This represents the fact that once all $4$ values are taken, the $132$ constraint is no longer relevent.
\end{enumerate}

An example of these operating in practice is given in figure~\ref{fig:all4}.

\tikzset{mydot/.style={circle,fill=black,draw=black,inner sep=1pt}}
\tikzset{mynode/.style={name=#1,draw=black,rounded corners,inner sep=2pt,align=center}}

\newcommand\perm[3]{
	\tikz{
		\begin{scope}[scale=0.17]
			\foreach[count=\i] \j in {#1} \node[mydot] at (\j,-\i) {};
			\foreach \j in {#2} \draw (\j+0.5,-0.5) -- (\j+0.5,-\i-0.5);
			\ifx&#1&%
			\def\i{1}
			\fi
			\draw (0.5,-0.5) -- (0.5,-\i-0.5);
			\begin{scope}[shift={(0.5+\i,-0.5*\i)},scale=4]
				\ifx&#3&%
				\node at (0.5,-0.1) {$\varnothing$};
				\else
				\linkpattern{#3}
				\fi
			\end{scope}
		\end{scope}
}}

\begin{figure}
	\tikzset{mylab/.style={fill=white,circle,draw=black,inner sep=0.8pt,align=center,node font=\small}}
	\begin{tikzpicture}[rotate=90,xscale=0.93,yscale=2.8]
	\node[mynode=a] at (0,0) {\perm{}{}{}};
	\node[mynode=b] at (-5,-1) {\perm{1}{}{}};
	\node[mynode=c] at (2,-1) {\perm{1}{1}{1/2}};
	\draw[-latex] (a) -- node[mylab] {$3$} (b);
	\draw[-latex] (a) -- node[mylab] {$1$} (c);
	\node[mynode=d1] at (-9,-2) {\perm{2,1}{}{}};
	\node[mynode=e1] at (-6.5,-2) {\perm{2,1}{1}{1/2}};
	\node[mynode=d2] at (-4,-2) {\perm{1,2}{}{}};
	\node[mynode=e2] at (5,-2) {\perm{1,2}{2}{1/2}};
	\node[mynode=e3] at (2,-2) {\perm{1,2}{1}{1/2}};
	\node[mynode=e4] at (-1,-2) {\perm{2,1}{2}{1/2}};
	\node[mynode=g1] at (7.4,-2) {\perm{2,1}{1,2}{1/2,3/4}};
	\node[mynode=g2] at (10.5,-2) {\perm{1,2}{1,2}{1/4,2/3}};
	\draw[-latex] (b) -- node[mylab] {$3$} (d1);
	\draw[-latex] (b) -- node[mylab] {$1$} (e1);
	\draw[-latex] (c) -- node[mylab] {$4$} (d2);
	\draw[-latex] (c) -- node[mylab] {$1$} (g1);
	\draw[-latex] (c) -- node[mylab] {$1$} (g2);
	\draw[-latex] (c) -- node[mylab] {$2$} (e2);
	\draw[-latex] (c) -- node[mylab] {$3$} (e3);
	\draw[-latex] (c) -- node[mylab] {$3$} (e4);
	\node[mynode=f1] at (-11,-3) {\perm{3,2,1}{}{}};
	\node[mynode=f1b] at (-10,-3) {\perm{3,2,1}{1}{1/2}};
	\node[mynode=f3] at (-9,-3) {\perm{3,1,2}{}{}};
	\node[mynode=f3b] at (-6,-3) {\perm{3,1,2}{2}{1/2}};
	\node[mynode=f3c] at (-7,-3) {\perm{3,1,2}{1}{1/2}};
	\node[mynode=f3d] at (-8,-3) {\perm{3,2,1}{2}{1/2}};
	\node[mynode=f2] at (-5,-3) {\perm{2,3,1}{}{}};
	\node[mynode=h1] at (-4,-3) {\perm{2,3,1}{1}{1/2}};
	\node[mynode=f4] at (5,-3) {\perm{1,2,3}{}{}};
	\node[mynode=f4b] at (8,-3) {\perm{1,2,3}{3}{1/2}};
	\node[mynode=f4c] at (7,-3) {\perm{1,2,3}{2}{1/2}};
	\node[mynode=f4d] at (6,-3) {\perm{2,3,1}{3}{1/2}};
	\node[mynode=f5] at (1,-3) {\perm{1,3,2}{}{}};
	\node[mynode=f5b] at (4,-3) {\perm{1,3,2}{2}{1/2}};
	\node[mynode=f5c] at (3,-3) {\perm{1,3,2}{1}{1/2}};
	\node[mynode=f5d] at (2,-3) {\perm{2,3,1}{2}{1/2}};
	\node[mynode=f6] at (-3,-3) {\perm{2,1,3}{}{}};
	\node[mynode=f6b] at (0,-3) {\perm{2,1,3}{3}{1/2}};
	\node[mynode=f6c] at (-1,-3) {\perm{2,1,3}{2}{1/2}};
	\node[mynode=f6d] at (-2,-3) {\perm{3,2,1}{3}{1/2}};
	\node[mynode=h2] at (9,-3) {\perm{3,1,2}{3}{1/2}};
	\node[mynode=h3] at (10,-3) {\perm{2,1,3}{1}{1/2}};
	\node[mynode=h4] at (11,-3) {\perm{1,2,3}{1}{1/2}};
	\foreach \x/\y/\z in {d1/f1/3,d1/f1b/1,d2/f2/3,e1/f3/4,e1/f3b/2,e1/f3c/3,e1/f3d/3,d2/h1/1,g1/h2/4,g1/h3/4,g2/h4/4,e2/f4/4,e2/f4b/2,e2/f4c/3,e2/f4d/3,e3/f5/4,e3/f5b/2,e3/f5c/3,e3/f5d/3,e4/f6/4,e4/f6b/2,e4/f6c/3,e4/f6d/3} \draw[-latex] (\x) -- node[mylab] {$\z$} (\y);
	\node[mynode=g1] at (-11,-4) {\perm{4,3,2,1}{}{}};
	\node[mynode=g1b] at (-10,-4) {\perm{4,3,1,2}{}{}};
	\node[mynode=g3] at (-9,-4) {\perm{4,2,3,1}{}{}};
	\node[mynode=g3b] at (-6,-4) {\perm{4,1,2,3}{}{}};
	\node[mynode=g3c] at (-7,-4) {\perm{4,1,3,2}{}{}};
	\node[mynode=g3d] at (-8,-4) {\perm{4,2,1,3}{}{}};
	\node[mynode=g2] at (-5,-4) {\perm{3,4,2,1}{}{}};
	\node[mynode=i1] at (-4,-4) {\perm{3,4,1,2}{}{}};
	\node[mynode=g4] at (5,-4) {\perm{2,3,4,1}{}{}};
	\node[mynode=g4b] at (8,-4) {\perm{1,2,3,4}{}{}};
	\node[mynode=g4c] at (7,-4) {\perm{1,2,4,3}{}{}};
	\node[mynode=g4d] at (6,-4) {\perm{2,3,1,4}{}{}};
	\node[mynode=g5] at (1,-4) {\perm{2,4,3,1}{}{}};
	\node[mynode=g5b] at (4,-4) {\perm{1,4,2,3}{}{}};
	\node[mynode=g5c] at (3,-4) {\perm{1,4,3,2}{}{}};
	\node[mynode=g5d] at (2,-4) {\perm{2,4,1,3}{}{}};
	\node[mynode=g6] at (-3,-4) {\perm{3,2,4,1}{}{}};
	\node[mynode=g6b] at (0,-4) {\perm{2,1,3,4}{}{}};
	\node[mynode=g6c] at (-1,-4) {\perm{2,1,4,3}{}{}};
	\node[mynode=g6d] at (-2,-4) {\perm{3,2,1,4}{}{}};
	\node[mynode=i2] at (9,-4) {\perm{3,1,2,4}{}{}};
	\node[mynode=i3] at (10,-4) {\perm{3,1,4,2}{}{}};
	\node[mynode=i4] at (11,-4) {\perm{1,3,4,2}{}{}};
	\foreach \x/\y/\z in {f1/g1/3,f1b/g1b,f3/g3/3,f3b/g3b,f3c/g3c,f3d/g3d,f2/g2/3,h1/i1,f4/g4/3,f4b/g4b,f4c/g4c,f4d/g4d,f5/g5/3,f5b/g5b,f5c/g5c,f5d/g5d,f6/g6/3,f6b/g6b,f6c/g6c,f6d/g6d,h2/i2,h3/i3,h4/i4} \draw[-latex] (\x) -- node[mylab]{$\ifx\y\z4\else3\fi$} (\y);
	\end{tikzpicture}
	\caption{Construction of all twenty three $Av(1324)$ PAPs of length $4$. The numbers in lines match the actions from section~\ref{sec:state}. The dot diagrams are the permutations being built up (read from top to bottom), with vertical bars indicating where future dots should be inserted; the link diagram is the state. The permutation diagrams should be read from the top down; the top right permutation diagram represents the permutation $1342$. The permutation $1324$ is avoided due to not allowing a second type $4$ transition from the top most state in the third column. }
	\label{fig:all4}
\end{figure}
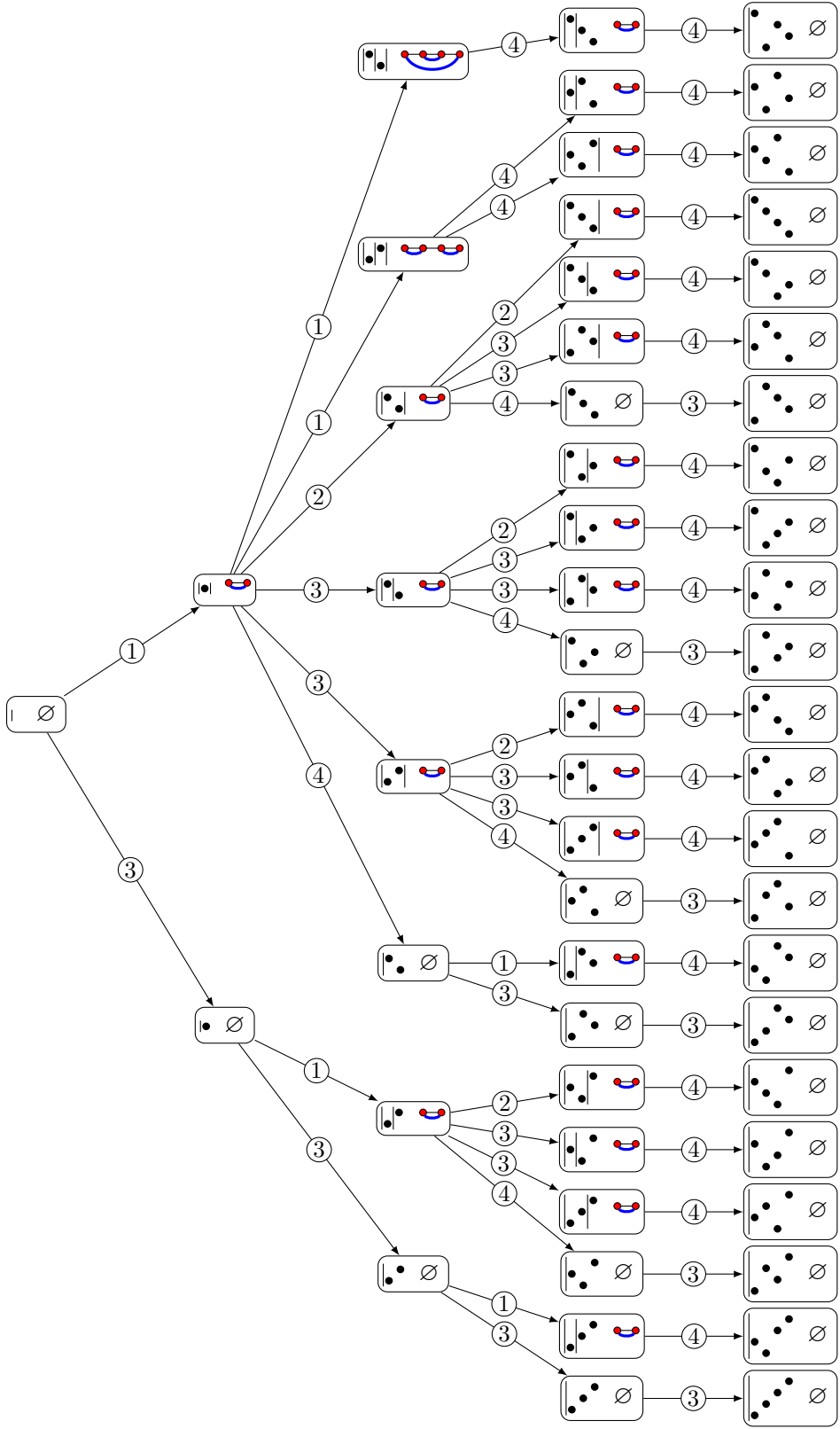

\newcommand\slp[1]{
	\tikz{
		\begin{scope}[scale=0.5]
			\ifx&#1&%
			\node at (0.5,-0.1) {$\varnothing$};
			\else
			\linkpattern{#1}
			\fi
		\end{scope}
}}

\newcommand\slpp[2]{
	\tikz{
		\begin{scope}[scale=0.5]
			\ifx&#1&%
			\node at (0.5,-0.1) {$\varnothing$};
			\else
			\linkpattern{#1}
			\fi
		\end{scope}
		\node at (0.25,-0.45) {#2}
}}

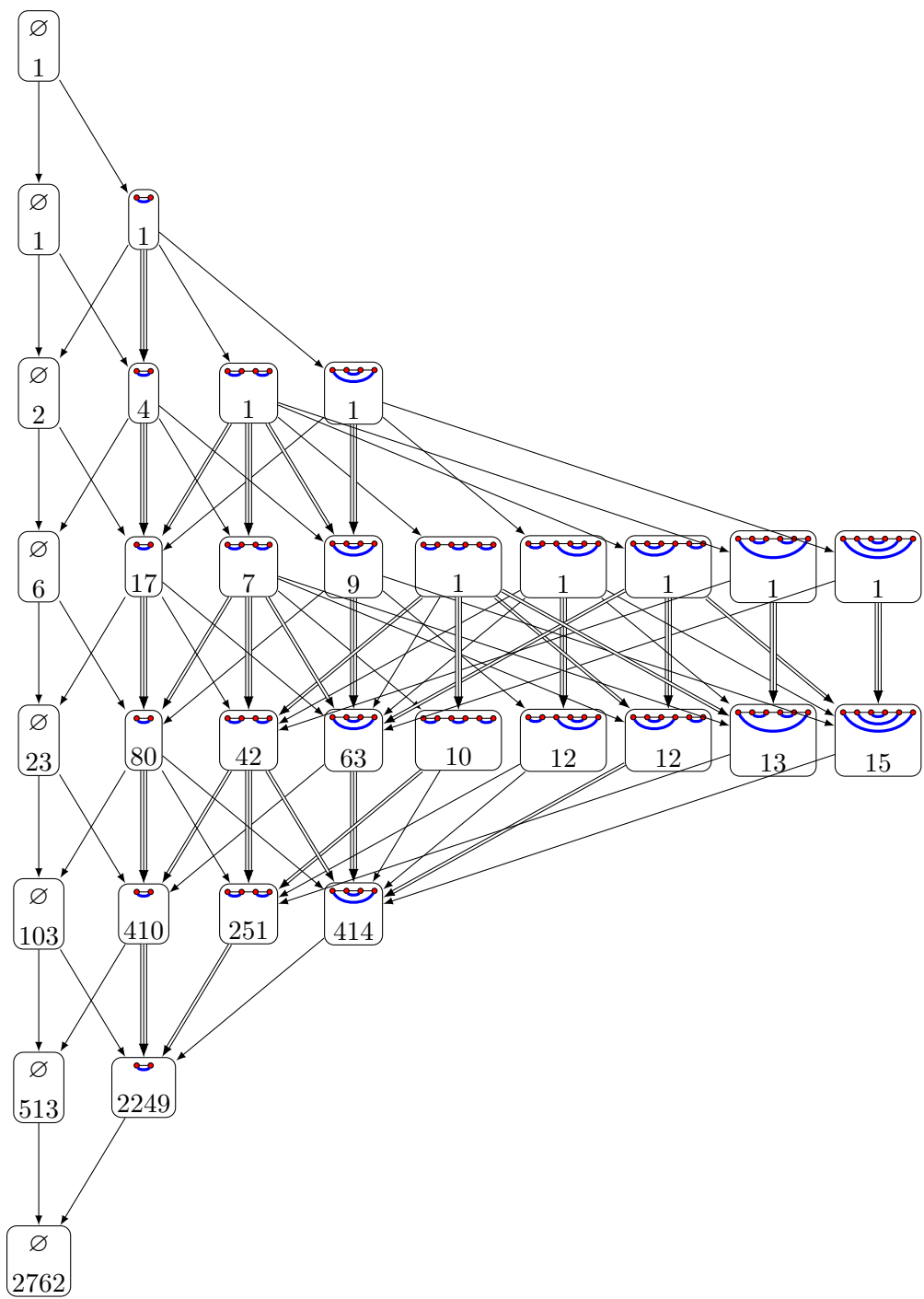
\begin{figure}
	\tikzset{mylab/.style={fill=white,circle,draw=black,inner sep=0.8pt}}
	\begin{tikzpicture}[xscale=1.5,yscale=2.5,double distance=0.8pt]
        \tikzset{triple/.style={preaction={double,double distance=2pt,draw}}}
	\node[mynode] (n8_) at (0,8) {\slp{} \\ $1$};
	\node[mynode] (n7_) at (0,7) {\slp{} \\ $1$};
	\node[mynode] (n7_0t1) at (1,7) {\slp{1/2} \\ $1$};
	\draw[-latex] (n8_) to (n7_0t1);
	\draw[-latex] (n8_) to (n7_);
	\node[mynode] (n6_) at (0,6) {\slp{} \\ $2$};
	\node[mynode] (n6_0t1) at (1,6) {\slp{1/2} \\ $4$};
	\node[mynode] (n6_0t1_2t3) at (2,6) {\slp{1/2,3/4} \\ $1$};
	\node[mynode] (n6_1t2_0t3) at (3,6) {\slp{2/3,1/4} \\ $1$};
	\draw[-latex] (n7_) to (n6_0t1);
	\draw[-latex,triple] (n7_0t1) to (n6_0t1);
	\draw[-latex] (n7_0t1) to (n6_);
	\draw[-latex] (n7_) to (n6_);
	\draw[-latex] (n7_0t1) to (n6_0t1_2t3);
	\draw[-latex] (n7_0t1) to (n6_1t2_0t3);
	\node[mynode] (n5_) at (0,5) {\slp{} \\ $6$};
	\node[mynode] (n5_0t1) at (1,5) {\slp{1/2} \\ $17$};
	\node[mynode] (n5_0t1_2t3) at (2,5) {\slp{1/2,3/4} \\ $7$};
	\node[mynode] (n5_1t2_0t3) at (3,5) {\slp{2/3,1/4} \\ $9$};
	\node[mynode] (n5_0t1_2t3_4t5) at (4,5) {\slp{1/2,3/4,5/6} \\ $1$};
	\node[mynode] (n5_0t1_3t4_2t5) at (5,5) {\slp{1/2,4/5,3/6} \\ $1$};
	\node[mynode] (n5_1t2_0t3_4t5) at (6,5) {\slp{2/3,1/4,5/6} \\ $1$};
	\node[mynode] (n5_1t2_3t4_0t5) at (7,5) {\slp{2/3,4/5,1/6} \\ $1$};
	\node[mynode] (n5_2t3_1t4_0t5) at (8,5) {\slp{3/4,2/5,1/6} \\ $1$};
	\draw[-latex] (n6_0t1_2t3) to (n5_0t1_2t3_4t5);
	\draw[-latex] (n6_) to (n5_0t1);
	\draw[-latex,triple] (n6_0t1) to (n5_0t1);
	\draw[-latex,double] (n6_0t1_2t3) to (n5_0t1);
	\draw[-latex,triple] (n6_1t2_0t3) to (n5_1t2_0t3);
	\draw[-latex] (n6_0t1_2t3) to (n5_1t2_0t3_4t5);
	\draw[-latex] (n6_0t1) to (n5_);
	\draw[-latex] (n6_1t2_0t3) to (n5_0t1);
	\draw[-latex] (n6_0t1_2t3) to (n5_1t2_3t4_0t5);
	\draw[-latex] (n6_1t2_0t3) to (n5_2t3_1t4_0t5);
	\draw[-latex] (n6_) to (n5_);
	\draw[-latex,triple] (n6_0t1_2t3) to (n5_0t1_2t3);
	\draw[-latex] (n6_1t2_0t3) to (n5_0t1_3t4_2t5);
	\draw[-latex,double] (n6_0t1_2t3) to (n5_1t2_0t3);
	\draw[-latex] (n6_0t1) to (n5_0t1_2t3);
	\draw[-latex] (n6_0t1) to (n5_1t2_0t3);
	\node[mynode] (n4_) at (0,4) {\slp{} \\ $23$};
	\node[mynode] (n4_0t1) at (1,4) {\slp{1/2} \\ $80$};
	\node[mynode] (n4_0t1_2t3) at (2,4) {\slp{1/2,3/4} \\ $42$};
	\node[mynode] (n4_1t2_0t3) at (3,4) {\slp{2/3,1/4} \\ $63$};
	\node[mynode] (n4_0t1_2t3_4t5) at (4,4) {\slp{1/2,3/4,5/6} \\ $10$};
	\node[mynode] (n4_0t1_3t4_2t5) at (5,4) {\slp{1/2,4/5,3/6} \\ $12$};
	\node[mynode] (n4_1t2_0t3_4t5) at (6,4) {\slp{2/3,1/4,5/6} \\ $12$};
	\node[mynode] (n4_1t2_3t4_0t5) at (7,4) {\slp{2/3,4/5,1/6} \\ $13$};
	\node[mynode] (n4_2t3_1t4_0t5) at (8,4) {\slp{3/4,2/5,1/6} \\ $15$};
	\draw[-latex] (n5_0t1_2t3) to (n4_0t1_2t3_4t5);
	\draw[-latex] (n5_2t3_1t4_0t5) to (n4_1t2_0t3);
	\draw[-latex] (n5_) to (n4_0t1);
	\draw[-latex,triple] (n5_0t1_3t4_2t5) to (n4_0t1_3t4_2t5);
	\draw[-latex,triple] (n5_1t2_0t3_4t5) to (n4_1t2_0t3_4t5);
	\draw[-latex] (n5_0t1_2t3_4t5) to (n4_1t2_0t3);
	\draw[-latex,triple] (n5_0t1) to (n4_0t1);
	\draw[-latex,double] (n5_0t1_2t3) to (n4_0t1);
	\draw[-latex] (n5_0t1_3t4_2t5) to (n4_1t2_3t4_0t5);
	\draw[-latex,double] (n5_0t1_2t3_4t5) to (n4_1t2_0t3_4t5);
	\draw[-latex,triple] (n5_1t2_0t3) to (n4_1t2_0t3);
	\draw[-latex,double] (n5_0t1_2t3_4t5) to (n4_0t1_2t3);
	\draw[-latex,triple] (n5_2t3_1t4_0t5) to (n4_2t3_1t4_0t5);
	\draw[-latex] (n5_0t1_2t3) to (n4_1t2_0t3_4t5);
	\draw[-latex] (n5_0t1) to (n4_);
	\draw[-latex] (n5_0t1_3t4_2t5) to (n4_2t3_1t4_0t5);
	\draw[-latex,double] (n5_1t2_0t3_4t5) to (n4_2t3_1t4_0t5);
	\draw[-latex,triple] (n5_1t2_3t4_0t5) to (n4_1t2_3t4_0t5);
	\draw[-latex] (n5_1t2_0t3) to (n4_0t1);
	\draw[-latex,double] (n5_1t2_0t3_4t5) to (n4_1t2_0t3);
	\draw[-latex] (n5_0t1_3t4_2t5) to (n4_0t1_2t3);
	\draw[-latex] (n5_0t1_3t4_2t5) to (n4_1t2_0t3);
	\draw[-latex] (n5_0t1_2t3) to (n4_1t2_3t4_0t5);
	\draw[-latex] (n5_1t2_0t3) to (n4_2t3_1t4_0t5);
	\draw[-latex] (n5_1t2_3t4_0t5) to (n4_0t1_2t3);
	\draw[-latex] (n5_) to (n4_);
	\draw[-latex,triple] (n5_0t1_2t3) to (n4_0t1_2t3);
	\draw[-latex] (n5_1t2_0t3) to (n4_0t1_3t4_2t5);
	\draw[-latex,double] (n5_0t1_2t3) to (n4_1t2_0t3);
	\draw[-latex] (n5_0t1) to (n4_0t1_2t3);
	\draw[-latex,double] (n5_0t1_2t3_4t5) to (n4_1t2_3t4_0t5);
	\draw[-latex,triple] (n5_0t1_2t3_4t5) to (n4_0t1_2t3_4t5);
	\draw[-latex] (n5_0t1) to (n4_1t2_0t3);
	\node[mynode] (n3_) at (0,3) {\slp{} \\ $103$};
	\node[mynode] (n3_0t1) at (1,3) {\slp{1/2} \\ $410$};
	\node[mynode] (n3_0t1_2t3) at (2,3) {\slp{1/2,3/4} \\ $251$};
	\node[mynode] (n3_1t2_0t3) at (3,3) {\slp{2/3,1/4} \\ $414$};
	\draw[-latex] (n4_2t3_1t4_0t5) to (n3_1t2_0t3);
	\draw[-latex] (n4_) to (n3_0t1);
	\draw[-latex] (n4_0t1_2t3_4t5) to (n3_1t2_0t3);
	\draw[-latex,triple] (n4_0t1) to (n3_0t1);
	\draw[-latex,double] (n4_0t1_2t3) to (n3_0t1);
	\draw[-latex,triple] (n4_1t2_0t3) to (n3_1t2_0t3);
	\draw[-latex,double] (n4_0t1_2t3_4t5) to (n3_0t1_2t3);
	\draw[-latex] (n4_0t1) to (n3_);
	\draw[-latex] (n4_1t2_0t3) to (n3_0t1);
	\draw[-latex,double] (n4_1t2_0t3_4t5) to (n3_1t2_0t3);
	\draw[-latex] (n4_0t1_3t4_2t5) to (n3_0t1_2t3);
	\draw[-latex] (n4_0t1_3t4_2t5) to (n3_1t2_0t3);
	\draw[-latex] (n4_1t2_3t4_0t5) to (n3_0t1_2t3);
	\draw[-latex] (n4_) to (n3_);
	\draw[-latex,triple] (n4_0t1_2t3) to (n3_0t1_2t3);
	\draw[-latex,double] (n4_0t1_2t3) to (n3_1t2_0t3);
	\draw[-latex] (n4_0t1) to (n3_0t1_2t3);
	\draw[-latex] (n4_0t1) to (n3_1t2_0t3);
	\node[mynode] (n2_) at (0,2) {\slp{} \\ $513$};
	\node[mynode] (n2_0t1) at (1,2) {\slp{1/2} \\ $2249$};
	\draw[-latex] (n3_) to (n2_0t1);
	\draw[-latex,triple] (n3_0t1) to (n2_0t1);
	\draw[-latex,double] (n3_0t1_2t3) to (n2_0t1);
	\draw[-latex] (n3_0t1) to (n2_);
	\draw[-latex] (n3_1t2_0t3) to (n2_0t1);
	\draw[-latex] (n3_) to (n2_);
	\node[mynode] (n1_) at (0,1) {\slp{} \\ $2762$};
	\draw[-latex] (n2_0t1) to (n1_);
	\draw[-latex] (n2_) to (n1_);
	\end{tikzpicture}
	\caption{The transfer matrix method for $Av(1324)$ PAPs of length $7$. Double and triple lines represent connections with weights of two and three respectively. }
	\label{fig:sum7}
\end{figure}


Note that a state with more links than remaining elements to insert can never be completed and should be discarded.
Inversely, each step can only increase the size of a link pattern by $2$.
In other words, we have the inequalities:
\begin{equation}\label{ineq}
k\le \min(s,t)
\qquad
\text{where }
\begin{cases}
k=\text{half-size of link patterns},
\\
s=\#\text{elements left to insert},
\\
t=\#\text{past time steps}
\end{cases}
\end{equation}
and we recall that $s+t=n$.

A link diagram state can be mapped to a state similar to the states used in  \cite{CG15} by replacing each link by a bracket containing anything under said link, followed by a placeholder for one or more numbers. The whole thing should be preceded by a placeholder for 0 or more numbers.

\begin{remark*}
If one uses only moves (1) and (4), one generates this way (the number
of) {\em alternating}\/ $1324$-avoiding permutations \cite{A217807}.
\end{remark*}

\subsection{Transfer matrix}

The set of equations defined by the state and child states above can be solved recursively in a straightforward manner using dynamic programming, but a
more memory efficient method is to use the transfer matrix method. This involves two sets at each point, a source and
a destination. Start with the ``source'' set consisting of the start state $\varnothing$ and associated multiplicity 1. Then
for each child state of each element in the source set, add that to a ``destination'' set with the same multiplicity as associated
with the source state. When an element gets added multiple times to the destination set, sum the multiplicities. 
Each element in the destination set will now have $n-1$ elements remaining. Now clear the (no longer needed) source
set, swap the source and destination set, and repeat the process until down to $0$ remaining elements. The multiplicity of the
$\varnothing$ element is the desired answer.

An example of the full transfer matrix computation for length seven
is given in figure~\ref{fig:sum7}. Each row represents one set
in the transfer matrix. Each node contains the state and the associated multiplicity.
Lines between them indicate that the lower state is a child of the upper state; if it is a child multiple ways,
that multiplicity is shown through the thickness of the line. The number of permutations of $n$ elements can
be found from the multiplicity of the $\varnothing$ node on line $n+1$ of this figure.

The transfer matrix saves some memory relative to dynamic programming as old states can be discarded when not needed. 
Since each state in a set contains exactly the
same number of remaining elements, there is no need to store that as part of the state. As each link diagram can be mapped 
to an integer efficiently, there is no need to use a hash table, rather a simple array can be used.

The main algorithm then consists of passing elements from one array to another. Instead of iterating over the source
array and computing all destinations, it is possible to do the converse; to iterate over the destination array
and compute all sources. This is mathematically equivalent, but in the case of parallel processing in a shared
memory system it is significantly more efficient as each thread can be assigned some destination nodes to process,
and may write to said nodes without computationally expensive synchronization with other threads for each write.
This allows almost perfect efficiency of multi-threading on a single shared memory computer. However, this approach
is counter-productive in a message passing parallelisation system, where a shared write is fast (send a message to
the node storing that state), while a shared read is slow (send a message to the node storing that state, and wait for a response or 
go on to do something else, saving enough information to be able to remember what to do with the result when it does return).

The performance (both time and memory) is basically determined by the number of states, which is the sum of the Catalan numbers up to size $n/2$ (rounded down).
This is obtained by maximizing $k$ in Eq.~\eqref{ineq}, i.e., after $n/2$ iterations the longest states will be those with $n/2$ links\footnote{These come from adding $1$, $3$, $5$, up to $n-1$, in some 132 avoiding permutation. This means there will be exactly one multiplicity for each state in the first step of the transfer matrix containing that state. Given this, there is no point storing or calculating these states, reducing the memory use by effectively one term. We used this, and also the next term -- states which had one adjacency -- to get yet another term without extra memory use. The formula for that is somewhat more complex.}. Afterwards, the number of states will reduce as they will
be constrained by the number of remaining elements, being an upper bound on the number of links allowed. As the Catalan numbers grow like $4^n/\sqrt{\pi n^3}$, the
algorithm uses roughly $4$ times as much memory each time $n$ increases by $2$.

\subsection{Running}

We wrote a C program using message passing to run on a distributed system, and ran it in the Spartan \cite{Spartan} cluster at the High Performance Computing Centre at the University of Melbourne on 168 cores with 20GB per core. The computation was performed five times, each with computations
performed modulo a number close to $2^{16}$, so that only $16$ bits of storage were needed for each state. Each run took several hours. The actual answers were then reconstituted using the Chinese remainder theorem. This produced the series up to length $50$ permutations. The series is presented
in table~\ref{table:series} (see also \cite{A061552}).

\begin{table}
{\renewcommand{\arraystretch}{0.95}
	\begin{tabular}{l}
1  \\
2  \\
6  \\
23  \\
103  \\
513  \\
2762  \\
15793  \\
94776  \\
591950  \\
3824112  \\
25431452  \\
173453058  \\
1209639642  \\
8604450011  \\
62300851632  \\
458374397312  \\
3421888118907  \\
25887131596018  \\
198244731603623  \\
1535346218316422  \\
12015325816028313  \\
94944352095728825  \\
757046484552152932  \\
6087537591051072864  \\
49339914891701589053  \\
402890652358573525928  \\
3313004165660965754922  \\
27424185239545986820514  \\
228437994561962363104048  \\
1914189093351633702834757  \\
16130725510342551986540152  \\
136664757387536091240503406  \\
1163812341034817216384582333  \\
9959364766841851088593974979  \\
85626551244475524038311935717  \\
739479176041581588794042743521  \\
6413612398452364144369673970347  \\
55855094052029166019855630997080  \\
488354507551082299792086219184434  \\
4286013140398612535730177106798038  \\
37753338738386034300928290519149333  \\
333720028221302436110132711265898937  \\
2959914488410727889919188039470296624  \\
26338690757116988316771828238926079326  \\
235113956679181729949424482617740434207  \\
2105162587512716675745868833684827184388  \\
18904804517351837590874336467009693522354  \\
170253750251391700942449152528030601519757  \\
1537516984674177479234766336099763469212469 \\
	\end{tabular}
}	\caption{The first $50$ terms of the $Av(1324)$ series.}
	\label{table:series}
\end{table}

\section{Differential approximant analysis}
The most successful numerical method for extracting the asymptotics from the first few terms of the OGF of a function with a power-law singularity is the method of differential approximants, due to Guttmann and Joyce \cite{GJ72}, with subsequent refinements due to Baker and Hunter \cite{HB79} and Fisher and Au-Yang \cite{FA79}. Details are given in \cite{G14, GJ09, G89}. In brief, the method fits available coefficients to a judiciously chosen family of D-finite ordinary differential equations (ODEs), and the singularity structure of the ODEs is extracted by standard methods \cite{Ince27, Forsyth02}. 

For models with an isolated power-law singularity, the method is very successful, with estimates of the radius of convergence and critical exponents frequently accurate to $10$--$20$ significant digits or even more. However, when the method is used to analyse models with singularities that are not simple power-laws, such as those whose coefficients have stretched-exponential behaviour, the method fails, though in a predictable manner. In that case one finds that the radius of convergence estimates are typically only accurate to two or three significant digits, and the critical exponent estimates are numerically large and sporadic, and cannot be relied upon.

In this way, the method is useful -- as is a canary in a coal mine. If one analyses the known terms of a series with the method of differential approximants and finds estimates of the radius of convergence to be poorly converged, with numerically large exponent values, which are not mutually consistent, one can be confident that the underlying OGF does not have a power-law singularity. Applying the method to the first 30 terms of the $Av(1342)$ and $Av(1234)$ PAPs, which are known to be D-finite, the complete solution is found. Applying the method to the coefficients we have for $Av(1324)$ PAPs, the method suggests that the radius of convergence is around $0.09$, with an exponent variously estimated to be $-20$ or $+15$ or anything in between! This is the hallmark of a non-power-law singularity. 

However there is another property of differential approximants that has recently been elucidated that is very useful for our purposes. Trivially, every differential approximant is an ODE. That ODE reproduces exactly those coefficients used in its construction (obviously), but also implies {\em approximate}\/ values for all subsequent coefficients. What is surprising is how accurate these predicted coefficients can be {\em even when the differential approximants give poor estimates of the critical point and exponent}. We refer to this prediction of additional coefficients beyond those used in construction of the DA as {\em series extension}.

  The idea of series extension was introduced and studied in \cite{G16}, and has been used in several recent calculations. For example, in \cite{EPG17}, the authors took the first 32 coefficients of the cogrowth series of the group ${\mathbb Z} \wr {\mathbb Z}$ discussed above, and used them to predict a further 89 coefficients. The predicted error was obtained by constructing a number of differential approximants and, for each coefficient, forming the average among the predicted values. Taking the confidence-interval as 1.5 standard deviations from the mean, they showed that the actual error was less than or equal to the predicted error. The actual error increased from 1 part in $2.7 \times 10^{17}$ for the first predicted coefficient to 1 part in $4.6 \times 10^{3}$ for the eighty-ninth coefficient. In another example where the solution is not D-finite, that of 4-valent Eulerian orientations, counted by vertices\footnote{This is equivalent to the 6-vertex problem on a random quadrangulation.}, the authors used the first 100 coefficients to usefully predict the next 1000 with an accuracy of better than 1 part in $10^{28}$ \cite{EPG17b}.

The first use of this method here was to calculate the last coefficient. Using 16-bit co-primes, ten runs would have been needed to calculate $p_{50}.$ The method of series extension allowed us to predict the first {\em twenty-nine} digits of $p_{50},$ so only three primes were needed to calculate the coefficient, and we used a fourth as an extra check. For $Av(1324)$ PAPs we have 50 coefficients, and have used these to estimate the next 200 ratios and the next 175 coefficients. As discussed in \cite{G16}, the errors (estimated as 1.5 standard deviations) of the extrapolated ratios increase more slowly than those of the coefficients, so we can predict more estimated ratios than coefficients.

Further details of the method of differential approximants, its successes and limitations are discussed in \cite{G14}. Further details of its application to series extension  can be found in \cite{G16}.

\subsection{Ratio analysis}
As in our previous analysis of the 36-term series, our primary tool is based on the behaviour of the ratio of successive coefficients. We also make use of the approximate ratios and coefficients, as calculated by the method of series extension.

In the case of a simple power-law singularity with the asymptotic form of the coefficients given by  $a_n \sim const. \cdot \mu^n \cdot n^g$, the ratio of the coefficients is
\BE \label{eq:rat1}
r_n = \frac{a_n}{a_{n-1}} = \mu\left (1 + \frac{g}{n} +{\rm O}(\frac{1}{n^2}) \right ).
\EE
If on the other hand the coefficients  behave as 
\[
b_n \sim B \cdot \mu^n \cdot \mu_1^{n^\sigma} \cdot n^g,
\]
then the ratio of successive coefficients $r_n = b_n/b_{n-1}$, is
\begin{multline} \label{eq:rn}
r_n = \mu \left (1 + \frac{\sigma \log \mu_1}{n^{1-\sigma}} + \frac{g}{n} + \frac{\sigma^2 \log^2 \mu_1}{2n^{2-2\sigma}} + \frac {(\sigma-\sigma^2)\log \mu_1+2g\sigma \log \mu_1}{2n^{2-\sigma}} \right . \\
 \left . {}+ \frac{\sigma^3 \log^3 \mu_1}{6n^{3-3\sigma}} +{\rm O}(n^{2\sigma-3}) + {\rm O}(n^{-2}) \right ).
\end{multline}
In particular,
when $\sigma = \frac{1}{2}$, (here we are anticipating our findings), this specialises to
\BE \label{eq:half}
r_n = \mu \left (1 + \frac{ \log \mu_1}{2\sqrt{n}} + \frac{g+\frac{1}{8}\log^2 \mu_1}{n} + \frac{\log^3\mu_1+(6+24g)\log \mu_1 }{48n^{3/2} } + {\rm O}(n^{-2}) \right ).
\EE

\begin{figure}
\begin{minipage}{0.45\textwidth}
\centerline{\includegraphics[width=7cm]{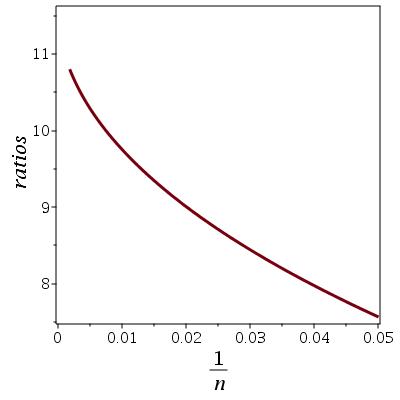} }
\caption{\label{fig:papr1}
Plot of ratios of coefficients against $\frac{1}{n}$.}
\end{minipage}
\begin{minipage}{0.45\textwidth}
\centerline{\includegraphics[width=7cm]{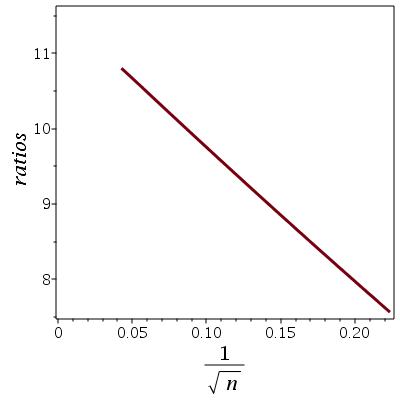} }
\caption{\label{fig:papr2}
Plot of ratios of coefficients  against $\frac{1}{\sqrt{n}}$.}
\end{minipage}
\end{figure}

In order to determine the nature of the asymptotic form of the coefficients of the $Av(1324)$  OGF, 
we first plot the ratios of successive coefficients $r_n = p_n/p_{n-1}$ against $1/n$, as shown in Figure~\ref{fig:papr1}.  In this and subsequent plots we have used the first 50 exact ratios and the subsequent 200 predicted ratios. Significant curvature is observed. This is inconsistent with an algebraic singularity, as can be seen from eqn.~\eqref{eq:rat1}. We next plot the same ratios against $1/\sqrt{n}$ in 
Figure~\ref{fig:papr2}, and the plot is seen to be visually linear, implying, from eqn.~\eqref{eq:rn} that $\sigma \approx 1/2.$  Linear extrapolation implies a limiting value as $n \to \infty$ around 11.60. 

\begin{figure}
\begin{minipage}{0.45\textwidth}
\centerline{\includegraphics[width=6.6cm]{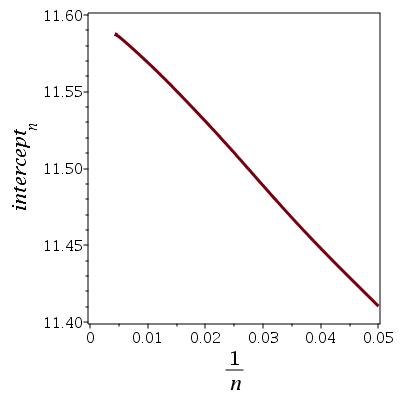} }
\caption{\label{fig:papextrap}
Linear extrapolants of successive ratios against $\frac{1}{n}$.}
\end{minipage}
\begin{minipage}{0.45\textwidth}
\centerline{\includegraphics[width=6.6cm]{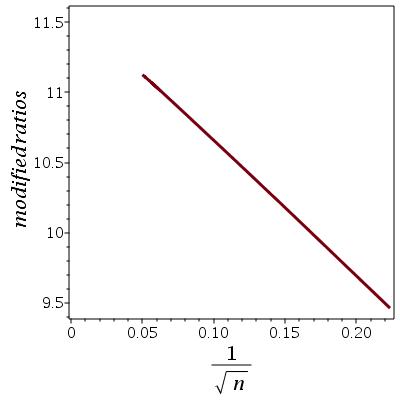} }
\caption{\label{fig:papll}
Modified ratios \eqref{eqn:linint} against $1/\sqrt{n}$.}
\end{minipage}
\end{figure}

We can significantly improve on this estimate by considering the sequence of extrapolants defined by successive pairs of points. That is to say, one can simply linearly extrapolate successive pairs of ratios $(r_k,r_{k+1})$ with $k$ increasing up to 240\footnote{Beyond this value, errors in the extrapolated ratios start to become visible.}. A plot of successive extrapolants against $1/n$ is shown in figure \ref{fig:papextrap}, which appears to be approaching a limit of around 11.60, or slightly below. We also take $\sigma =1/2$ as our (initial) conjectured value.

We can get somewhat more rapid convergence of the ratios if we remove the term $O(1/n)$ from eqn.~\eqref{eq:rn}, and this we do by calculating the modified ratios 
\begin{equation}\label{eqn:linint}
l_n = n \cdot r_n - (n-1)\cdot r_{n-1}= \mu \left (1+ \frac{\log {\mu_1}}{4\sqrt{n}} + O(n^{-3/2}) \right ).
\end{equation}

\begin{figure}
\begin{minipage}{0.45\textwidth}
\centerline{\includegraphics[width=6.6cm]{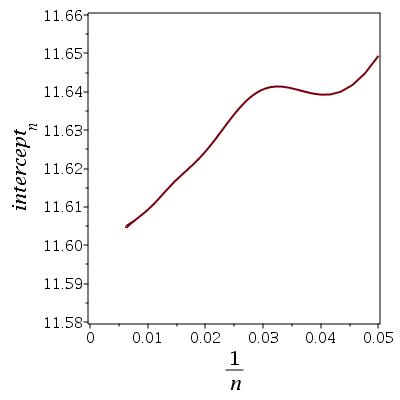} }
\caption{\label{fig:paplxt}
Estimates of  $\mu$ by linear extrapolants against $\frac{1}{n}$.}
\end{minipage}
\begin{minipage}{0.45\textwidth}
\centerline{\includegraphics[width=6.6cm]{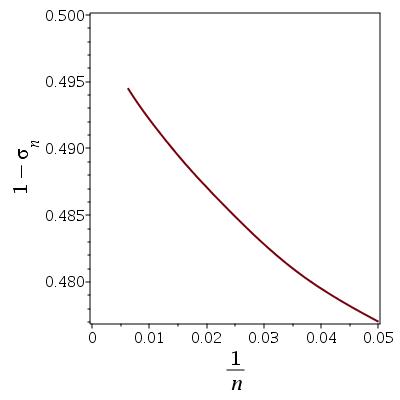} }
\caption{\label{fig:paps1}
Plot of estimates of  $1-\sigma$  against $\frac{1}{n}$.}
\end{minipage}
\end{figure}

We show the plot of these modified ratios $l_n$ against $1/\sqrt{n}$ in Figure~\ref{fig:papll}, and the linear extrapolants of successive pairs of points taken from Figure~\ref{fig:papll} (as described above) is shown in Figure~\ref{fig:paplxt}. We again conclude that $\mu$ is very close to, and probably slightly below, 11.60.

In order to more accurately estimate the value of the exponent $\sigma$, we note from (\ref{eq:rn}) that 
\[
(r_n/\mu-1) \sim const. n^{\sigma-1}.
\]
A log-log plot of $(1-r_n/\mu)$ against $\log{n},$ where we have taken $11.60$ as the value of $\mu$, is an uninteresting linear plot. However if we calculate the gradient from successive pairs of points, then the negative of this gradient is an estimator of the exponent $1-\sigma$.
We plot these estimators against $1/n$ in Figure~\ref{fig:paps1}, which provides compelling evidence that $\sigma = 1/2$. In our subsequent analysis, we will assume this value. Repeating this analysis with various values of $\mu$ around $11.60$, we find that a value slightly below this, around $\mu \approx 11.598$ is most consistent with $\sigma=1/2$.

Assuming then that $\sigma=1/2$, from \eqref{eq:half} it follows that 
\[
r_n/\mu = 1 + \frac{ \log \mu_1}{2\sqrt{n}} + \frac{g+\frac{1}{8}\log^2 \mu_1}{n} + {\rm O}(n^{-3/2}). 
\]

In order to estimate $\mu_1$ and $g$, we solve, sequentially, the equations 
\BE \label{eq:fit2a}
r_j/\mu = 1 + \frac{ c_1}{\sqrt{j}} + \frac{c_2}{j} +  \frac{c_3}{j^{3/2}}, 
\EE
for $j=k-1$, $j=k$ and $j=k+1$, with $k$ ranging from $3$ up to $49$.

\begin{figure}
\begin{minipage}{0.45\textwidth}
\centerline{\includegraphics[width=6.6cm]{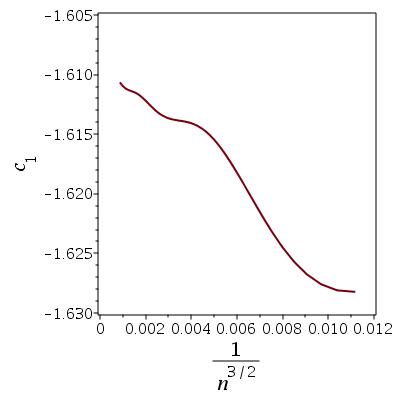} }
\caption{\label{fig:papc1}
Plot of estimates of parameter $c_1$ of (\ref{eq:fit2a}) against $\frac{1}{n^{3/2}}$.}
\end{minipage}
\begin{minipage}{0.45\textwidth}
\centerline{\includegraphics[width=6.6cm]{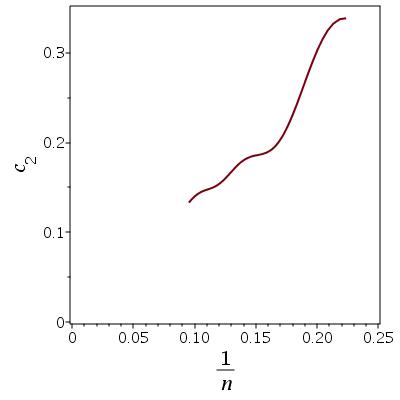} }
\caption{\label{fig:papc2}
Plot of estimates of parameter $c_2$ of (\ref{eq:fit2a}) against $\frac{1}{n}$.}
\end{minipage}
\end{figure}

The results are shown in figures \ref{fig:papc1} and \ref{fig:papc2}, plotting the parameters $c_1$ and $c_2$ respectively. The first neglected term in the asymptotics is O$(1/n^2)$ which is O$(1/n^{3/2})$ smaller than the term with coefficient $c_1,$ so $c_1$ is plotted against $1/n^{3/2}$. By a similar argument, $c_2$ is plotted against $1/{n}$. A simple visual extrapolation gives the estimate $c_1 = -1.6075 \pm 0.0025$. The plot for $c_2$ is difficult to extrapolate. It appears to be turning near its end point, and we very tentatively estimate  $0 < c_2 < 0.1$. From \eqref{eq:half}, $c_1=\log{\mu_1}/2$ and $c_2=g+\frac{1}{8}\log^2 \mu_1$. Hence we estimate  $\log{\mu_1} \approx -3.215,$ and $g \approx -1.2$. 

Next we tried direct fitting to the parameters in the assumed asymptotic form. That is to say, the assumed asymptotic form is 
\[
b_n \sim B\cdot \mu^n \cdot \mu_1^{n^\sigma} \cdot n^g.
\]
Therefore
\BE \label{logcan1}
\log {b_n} \sim \log{B} + n \log{\mu} + n^\sigma \log{\mu_1} + g \log {n}.
\EE
So if $\sigma$ is known, or assumed, we have four unknowns in this linear equation. It is then straightforward to solve the linear system
\[
\log {b_k} =  c_1k  + c_2 k^\sigma  +c_3 \log {k}+c_4
\]
for $k=n-2,\, n-1, \, n, \, n+1$ with $n$ ranging from $3$ up to the maximal number of known or approximated coefficients. Then $c_1$ estimates $\log {\mu}$, $c_2$ estimates $\log{\mu_1}$, $c_3$ estimates $g$ and $c_4$ gives estimators of $\log{B}$. 

An obvious useful variation is in those cases where, say, $\mu$ is known, or accurately estimated. Then one can solve  
\BE \label{eq:three}
\log(b_n) - n \log{\mu} =  c_2 n^\sigma  +c_3 \log {n} +c_4
\EE
 from three successive coefficients $b_{n-1}, \,\, b_n, \,\, b_{n+1}$, as before increasing the order of the lowest used coefficient by one until one runs out of coefficients. We do this below with our best estimate of $\mu.$

Fitting the available coefficients to the four unknowns, we estimate $c_1\approx 2.451$, implying $\mu\approx 11.60$, (consistent with our earlier estimate), $c_2=-3.225$, implying $\mu_1 \approx \exp(-3.225)=0.0398\ldots$,   $c_3\approx -1.15$, which is an estimate of the exponent $g$, and $c_4 \approx 1.7$ implying $B \approx 5.5$. 

If we set $\mu=11.598$, and fit to the remaining three unknowns, we find $c_2=-3.228$, implying $\mu_1 \approx 0.0396$,   $c_3\approx -1.15$, which is an estimate of the exponent $g$, and $c_4 \approx 1.9$ implying $B \approx 6.7$. 

This estimate of the power-law exponent $g$ around $-1.15$ needs further discussion. A pure power-law logarithmic term, such as $\log(1 - \mu \cdot z)$ would correspond to $g=-1$, as indeed would $z/\log(1 - \mu \cdot z)$. Power-law behaviour of the latter type (the reciprocal of a logarithm) is known to be difficult to analyse \cite{EPG17b}, typically giving rise to an exponent estimate of around $-1.2$ with most methods of analysis. For this reason we are very cautious in our estimate of the value of the exponent $g$, and so quote $g = -1.1 \pm 0.1$, to include the possibility that in fact $g=-1$ exactly. Indeed, an alternative analysis in Subsection~\ref{da} produces an estimate of $g$ which is closer to $-1$ than this analysis. That being said, if the estimate $g \approx -1.15$ can be believed, then seeking the simplest rational fraction would lead to the guesstimate $g = -7/6.$

\subsection{More complex stretched-exponential term?}
As noted above, the stretched exponential exponent can sometimes occur with an exponent that includes a multiplicative logarithmic factor, as seen in the cogrowth series for certain groups \cite{PS-C02}. With the extended series we have, we can test for this as follows:
Assume an asymptotic form for the coefficients of
\[
c_n \sim c\cdot \mu^n \cdot \mu_1^{\sqrt{n}\log^\delta{n}} \cdot n^g.
\]
Assume we know, or have a good estimate for, the growth constant $\mu$. Then form the normalised coefficients $d_n=c_n/\mu^n$. Then define
\[
e_n \equiv d_{n^2} \sim c \cdot \mu_2^{n\log^\delta{n}}\cdot n^{2g},
\]
where $\mu_2=\mu_1\cdot 2^\delta$. The ratio of successive coefficients is now
\[
{\tilde r}_n = \frac{e_n}{e_{n-1}}=\mu_2^{\log^\delta{n}}\cdot \left ( 1 + \frac{2g}{n} \right ).
\]

\begin{figure}
\centerline{\includegraphics[width=6.6cm]{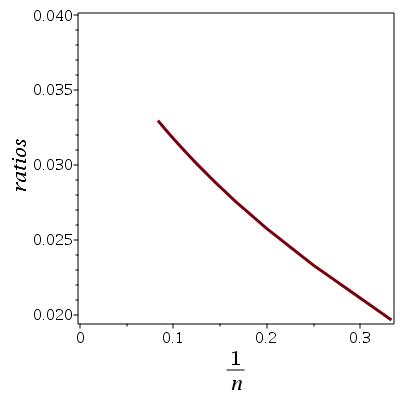} }
\caption{\label{fig:mu2}
Ratio plot of ${\tilde r}_n$ vs $1/n$, giving estimates of  $\mu_2$. }
\end{figure}

So if $\delta =0$, a plot of the ratios against $1/n$ should go to $\mu_1$ as $n \to \infty$; but if $\delta > 0$, the ratio plot should diverge; while if $\delta < 0$, the ratios should go to 1 as $n \to \infty$. Using the estimate $\mu=11.60$ the resultant ratio plot is shown in Figure~\ref{fig:mu2}, which is clearly extrapolating to a value around $0.04$, which is precisely the value previously estimated for $\mu_1$. Furthermore, this plot is not hugely sensitive to the estimate of $\mu$. Therefore we find that there is no evidence for a confluent logarithmic term in the exponent of the stretched-exponential.

\subsection{Using differential approximants}\label{da}

As noted above, differential approximants are useful insofar as they indicate that the singularity is not a pure power-law. They provide a signal, but are then of no further use in their current form. The presence of the $\mu_1^{\sqrt{n}}$ term is responsible for the lack of applicability of the method. However we can manipulate the series to remove the offending term, and then use this powerful method. From eqn.~\eqref{logcan1} one has

\[
\log {b_n} \sim \log{B} + n \log{\mu} + \sqrt{n} \log{\mu_1} + g \log {n}.
\]
Then defining $\tilde{b}_n = \log(b_n)/\sqrt{n}$, one has
\BE  \label{eq:renorm}
c_n=2n^{3/2}(\tilde{b}_n - \tilde{b}_{n-1}) = \left (2g-\log{B}+\frac{1}{4}\cdot\log{\mu}\right ) + n\log{\mu}-g\log{n}
\EE
So 
\[
d_n = \exp(c_n) =D\cdot \mu^n \cdot n^{-g}.
\]

In this way we have transformed the series to one whose coefficients behave asymptotically, at least to leading order, like a function with an algebraic singularity. We can therefore analyze the series with coefficients $d_n$ by the method of differential approximants (DAs). We summarise the results in table~\ref{tab:ana}. The column labelled $L$ gives the degree of the inhomogeneous polynomial of the approximating ODE. The entries give an average value for the position and exponent of the singularity of the fitted ODEs. Full details of the method are given in \cite{G89,GJ09}. It is seen that the 2nd-order DAs give estimates of the radius of convergence that can be summarised as $1/\mu = 0.086205 \pm 0.00003$, or $\mu = 11.600 \pm 0.004$. The 3rd-order approximants are more stable, allowing us to estimate $
\mu = 0.086206 \pm 0.00001$, or $\mu = 11.6001 \pm 0.0013$. This is remarkably close to estimates above, obtained by various other means. Note too that the estimate of the exponent $g$ is now closer to $-1$ than found by ratio methods, and, as discussed above, could possibly be $-1$ exactly.

\begin{table}
\caption{\label{tab:ana}
Critical point and exponent estimates for renormalised $Av(1324)$ PAPs}
\begin{center}
\begin{tabular}{lllll} \hline \hline
 $L$   &  \multicolumn{2}{c}{Second-order DA} &
       \multicolumn{2}{c}{Third-order DA} \\
\hline
    &  \multicolumn{1}{c}{$1/\mu$} & \multicolumn{1}{c}{$g-1$} &
      \multicolumn{1}{c}{$1/\mu$} & \multicolumn{1}{c}{$g-1$} \\
\hline
 0 & 0.086197& -2.045 &  0.086206& -2.077 \\
 1 & 0.086210& -2.61 &  0.086205& -2.078 \\
2 & 0.086215& -2.173 & 0.086204& -2.067 \\
3 & 0.086206& -2.110 & 0.086205& -2.074 \\
4 & 0.086203& -2.094 & 0.086209& -2.084 \\
 5 & 0.086198& -2.055 &  0.086207& -2.077 \\
 6& 0.086199& -2.058 &  0.086206& -2.077 \\
7 & 0.086195& -2.043& 0.086204& -2.071 \\
8 & 0.086209& -2.18 & 0.086206& -2.075 \\
9 & 0.086203& -2.21 & 0.086205& -2.072 \\
10 & 0.086167& -2.6 & 0.086203& -2.062 \\
\hline \hline
\end{tabular}
\end{center}
\end{table}

We can also apply other traditional techniques to the transformed series. The ratios of successive terms ($d_n$) of the transformed series when plotted against $1/n$ are now linear, but as this doesn't give us a better estimate of $\mu$ than $\mu \approx 11.60,$ we don't show the figure here. 

So in conclusion we suggest that the coefficients of $Av(1324)$ PAPs behave as
\[
B\cdot \mu^n \mu_1^{\sqrt{n}} \cdot n^g,
\]
where $\mu=11.600 \pm 0.003,$ $\mu_1 = 0.0400 \pm 0.0005$, $g = -1.1 \pm 0.1$ while the estimate of $B $ depends sensitively on the precise value of $\mu,\,\,\mu_1$ and $g,$ so that if we vary these quantities over their uncertainty range, the subsequent estimate of $B$ ranges from 1 to 13, so we don't quote an estimate.

For other length-$4$ PAPs, the growth constant $\mu$ is an integer. It is clearly not an integer in this case. If it is a simple algebraic number -- and we have no compelling reason why it should be -- then we note that $9+3\sqrt{3}/2 = 11.598\ldots$ is indistinguishable from our numerical estimate. As discussed above, the exponent $g$ could be exactly $-1$ (corresponding to some power of a logarithm) or as high as $-6/5$, but the closest simple rational is $-7/6$. Despite these remaining uncertainties, this reanalysis using a longer series, and with series extension, provides compelling arguments for the presence of a stretched-exponential term of the form $\mu_1^{\sqrt{n}}$ in the case of $Av(1324)$ PAPs.

If we take our most-favoured estimates, $\mu = 11.598$, $\mu_1=0.0398$, $g=-7/6$, then $B \approx 9.0$, and the $1000^{\rm th}$ coefficient is $3.7 \times 10^{1017}$. Note that in \cite{ES13}, Steingr\'imsson gives a Monte Carlo estimate of the ratio $p_{1001}/p_{1000} \approx 11.011$. From our central estimates of $\mu$, $\mu_1$ and $g$, we find for this ratio $11.009$, which is far better agreement than could reasonably be expected.

\section*{Acknowledgements}
We wish to thank the High Performance Computing Centre at The University of Melbourne for access to the Spartan cluster, on which the bulk of the calculations in this paper were performed. We also wish to thank Dr Iwan Jensen for help in computing the last coefficient. PZJ was supported by ARC grant FT150100232.

\gdef\MRshorten#1 #2MRend{#1}%
\gdef\MRfirsttwo#1#2{\if#1M%
MR\else MR#1#2\fi}
\def\MRfix#1{\MRshorten\MRfirsttwo#1 MRend}
\newcommand\MR[1]{\relax\ifhmode\unskip\spacefactor3000 \space\fi
\MRhref{\MRfix{#1}}{{\scriptsize \MRfix{#1}}}}
\newcommand{\MRhref}[2]{%
\href{http://www.ams.org/mathscinet-getitem?mr=#1}{#2}}
\bibliographystyle{amsplainhyper}
\bibliography{a1324}

\end{document}

\end{document}